\documentclass[11pt, a4paper]{article}

\usepackage[utf8]{inputenc}
\usepackage[english]{babel}
\usepackage{amscd,amssymb}
\usepackage{amsmath}
\usepackage{amsfonts}
\usepackage{amsthm}
\usepackage{fancybox}
\usepackage{epic,eepic}
\usepackage{amstext}
\usepackage{mathtools}
\usepackage{esint}
\usepackage{booktabs}
\usepackage[hide links]{hyperref}
\usepackage{comment}
\usepackage{mathtools}
\usepackage{pgfplots}
\usepackage{todonotes}
\usepackage{tabularx}
\usepackage{enumerate}
\usepackage{bm}
\usepackage{stmaryrd}
\usepackage[noabbrev, capitalise, nameinlink]{cleveref}
\usepackage[style=numeric,
            sorting=none,
            giveninits=true,
            natbib=true]{biblatex}
\AtEveryBibitem{%
  \clearfield{issn} 
  \clearfield{month}

  \ifentrytype{online}{}{
    \clearfield{url}
  }
}

\crefname{equation}{}{}
\crefname{assumption}{Assumption}{Assumptions}
\crefformat{equation}{\textup{#2(#1)#3}}

\newtheorem{theorem}{Theorem}
\newtheorem{lemma}[theorem]{Lemma}
\newtheorem{corollary}[theorem]{Corollary}
\newtheorem{remark}[theorem]{Remark}

\crefname{corollary}{Corollary}{Corollaries}
\crefname{lemma}{Lemma}{Lemmas}

\renewcommand{\u}{u}
\renewcommand{\v}{v}

\newcommand{\uh}{u}
\newcommand{\vh}{v}

\newcommand{\numax}{\nu_{\mathrm{max}}}
\newcommand{\numin}{\nu_{\mathrm{min}}}

\newcommand{\uhe}{\u^e}

\newcommand{\uhG}{\u^{\mathrm{G}}}

\renewcommand{\L}[1]{L^2({#1})}
\newcommand{\Winf}[1]{W^{1,\infty}(#1)}

\newcommand{\mesh}{\mathbb{T}_h}

\newcommand{\pdisc}[2]{\mathcal{P}^{\mathrm{disc}}_{#2}(#1, \mesh)}
\newcommand{\pcont}[2]{\mathcal{P}_{#2}(#1, \mesh)}

\newcommand{\interop}{I_h}

\newcommand{\scal}[3]{(#1,#2)_{{#3}}}

\newcommand{\ah}{B}
\newcommand{\sumdg}[1]{\llbracket {#1} \rrbracket}
\newcommand{\sumdgw}[1]{\llbracket {#1} \rrbracket_w}
\newcommand{\jumpdg}[1]{\sumdg{{#1} \normal}}
\newcommand{\scaldg}[3]{{\langle {#1}, {#2} \rangle_{#3}}}
\newcommand{\ahosj}{B_{\osi}}

\newcommand{\ahos}{B_{\omega^*}}
\newcommand{\ahosi}{B_{\omega^*}^i}
\newcommand{\aposi}{B_{p,\omega^*}^i}
\newcommand{\acosi}{B_{c,\omega^*}^i}
\newcommand{\ahosb}{B_{\omega^*}^\partial}

\newcommand{\facedgint}{\mathbb{F}_h^i}
\newcommand{\facedgbdry}{\mathbb{F}_h^\partial}
\newcommand{\dgstab}{\gamma_h^2}
\newcommand{\tmin}{T_{\mathrm{min}}}
\newcommand{\tmax}{T_{\mathrm{max}}}
\newcommand{\reconop}{R_h}

\newcommand{\face}{F}

\newcommand{\lnorm}[2]{\|#1\|_{\L{#2}}}

\newcommand{\norm}[2]{\|{#1}\|_{#2}}

\newcommand{\normal}{\boldsymbol{n}}

\renewcommand{\P}{P}

\newcommand{\ddindex}{j}
\newcommand{\om}{\omega}
\newcommand{\omi}{\omega_{\ddindex}}
\newcommand{\osi}{\omega^{\ast}_{\ddindex}}
\newcommand{\os}{\omega^{\ast}}

\newcommand{\numDom}{M}

\newcommand{\omcon}{\kappa}
\newcommand{\oscon}{\kappa^{\ast}}

\newcommand{\eval}{\lambda}
\newcommand{\evec}{\varphi}
\newcommand{\evalindex}{k}
\newcommand{\evalloc}[1][\evalindex]{\eval_{\ddindex,#1}}
\newcommand{\evecloc}[1][\evalindex]{\evec_{\ddindex,#1}}

\newcommand{\Ccoer}{\alpha_{B}}
\newcommand{\Ccont}{C_{B}}

\def\dx{\,\text{d}x}

\allowdisplaybreaks

\title{Multiscale Spectral Generalized Finite Element Methods for Discontinuous Galerkin Schemes}

\author{Christian Alber\footnotemark[1]
\and Lukas Holbach\footnotemark[2]}

\date{\today}

\makeatletter
\let\@fnsymbol\@arabic
\makeatother

\begin{document}

\maketitle

\footnotetext[1]{Institute for Mathematics, Heidelberg University, Germany (\texttt{c.alber@uni-heidelberg.de}).}
\footnotetext[2]{Interdisciplinary Center for Scientific Computing, Heidelberg University, Germany\\ \hspace*{.6cm}(\texttt{lukas.holbach@iwr.uni-heidelberg.de}).}

\begin{abstract}
We propose a multiscale spectral generalized finite element method (MS-GFEM) for discontinuous Galerkin (DG) discretizations. The method builds local approximations on overlapping subdomains as the sum of a local source solution and a correction from an optimal spectral coarse space, which is obtained from a generalized eigenproblem. The global solution is then assembled via a partition of unity. We prove nearly exponential decay of the approximation error for second-order elliptic problems with highly heterogeneous diffusion discretized by a weighted symmetric interior-penalty DG scheme.
\end{abstract}



\section{Motivation and Model Setup}
\label{sec:intro}

This article presents a multiscale spectral generalized finite element method (MS-GFEM) tailored to discontinuous Galerkin (DG) discretizations of partial differential equations (PDEs). Our work is motivated by two problem classes: heterogeneous Stokes flows, where DG discretizations improve mass conservation, and convection-dominated diffusion, where DG fluxes enhance stability.
Previous work established nearly exponential error decay for MS-GFEM in the Hilbert-space setting \cite{Ma2025} and for continuous Galerkin methods \cite{Ma2022error}. As a first step toward extending that analysis to DG discretizations, we prove nearly exponential decay of the MS-GFEM approximation error for second-order elliptic PDEs with highly heterogeneous diffusion discretized by a weighted symmetric interior-penalty discontinuous Galerkin method.

Let $\Omega \subset \mathbb R^d$, $d \in \{2,3\}$, be a polygonal Lipschitz domain.   
We consider the following second-order elliptic PDE with homogeneous Dirichlet boundary conditions:
\vspace{-.2cm}
\begin{equation}
    \label{eq:PDEstrong}
    \left\{
    \begin{aligned}
        - \operatorname{div} (\nu \nabla u) &= f && \quad \text{in } \Omega,\\
        u &= 0 && \quad \text{on } \partial \Omega,
    \end{aligned}\right.
\end{equation}
where $f \in L^2(\Omega)$ and $\nu\in L^\infty(\Omega)$. 
We assume that $\numin \leq \nu(\bm x) \leq \numax$ holds for 
almost all $\bm x \in \Omega$ for some constants $\numax$, $\numin > 0$.
The weak form of~\cref{eq:PDEstrong} is based on the Sobolev space 
$H^1_{0}(\Omega)$ of $H^1$-functions vanishing on the boundary $\partial \Omega$.
Given $f\in L^2(\Omega)$, 
the weak form of \cref{eq:PDEstrong} seeks a function $u \in H^1_0(\Omega)$ such that 
\vspace{-.2cm}
\begin{equation}
    \label{eq:PDEweak}
    \int_\Omega \nu\nabla u \cdot\nabla v\dx = \int_\Omega fv\dx \eqqcolon \mathcal{F}(v) \qquad \forall \, v \in H^1_0(\Omega).  
\end{equation} 
By the Poincaré-Friedrichs inequality and the uniform bound on $\nu$, we obtain that 
the bilinear form in \cref{eq:PDEweak} defines an 
inner product on~$H^1_0(\Omega)$. Since $\mathcal{F}$ is a bounded linear functional, the Riesz representation theorem yields well-posedness of \cref{eq:PDEweak}.
For any subdomain $S \subset \Omega$, we define $\mathcal{F}_S(v) \coloneqq \int_S fv\dx$ for all $v \in L^2(S)$. We write $(\cdot, \cdot)_{S}$ and $\scaldg{\cdot}{\cdot}{F}$ for the $L^2$-inner products on subdomains $S$ and faces $F$, and use $\lesssim$ for inequalities up to constants independent of model parameters and mesh size.


\section{Discontinuous Galerkin Method}

For ease of presentation, we employ a matching shape-regular simplicial mesh $\mesh$
(cf. \cite[Section 1.4.1]{di2011mathematical}) and define $V_h \coloneqq \pdisc{\Omega}{1}$. 
For all $T\in\mesh$, we denote by $h_T$ the diameter of $T$ and define the mesh size $h \coloneqq \max_{T\in\mesh} h_T$.
We assume that the coefficient $\nu$ is resolved by the mesh, i.e., $\nu$ is piecewise constant with respect to $\mesh$.
The sets of interior and boundary faces of $\mesh$ are denoted by $\facedgint$ and $\facedgbdry$, respectively.
For an interior face $F = \partial T_1\cap \partial T_2$, we introduce $\nu_i = \nu|_{T_i}$, and we 
set $\nu_1=\nu_2 = \nu$ for boundary edges. 
To deal with varying $\nu$, we introduce $\nu$-dependent weights in the penalty and consistency terms (\cite[Chapter 4.5]{di2011mathematical}) and define the sum operators
\begin{equation*}
    \sumdgw{u} \coloneqq \frac{2\nu_2}{\nu_1 + \nu_2}u_1 + \frac{2\nu_1}{\nu_1 + \nu_2}u_2, 
    \qquad 
    \sumdg{u} \coloneqq u_1 + u_2, 
    \qquad u \in V_h^d.
\end{equation*}

For any subdomain $D\subset\Omega$, we define $\facedgint(D) \coloneqq \{ F\cap D : F\in \facedgint\}$ and
    $\facedgbdry(D) \coloneqq \{ F\cap \partial D : F\in \facedgbdry\}$, and we introduce the bilinear forms
\begin{equation*}
\begin{array}{l}
    B_{D}^i(\u, \v)=B_{p,D}^i(\u, \v)-B_{c,D}^i(\u, \v)-B_{c,D}^i(\v, \u), \\
    B_{D}^\partial(\u, \v)=B_{p,D}^\partial(\u, \v)-B_{c,D}^\partial(\u, \v)-B_{c,D}^\partial(\v, \u), \\
    B_{D}(\u, \v)=(\nu\nabla_h \u, \nabla_h \v)_{D}
    +B_{D}^i(\u, \v)+B_{D}^\partial(\u, \v),
\end{array}
\end{equation*}
with 
\begin{align*}
    &B_{c,D}^i(\u, \v)=\frac{1}{2} \scaldg{\sumdgw{\nu\nabla_h \u}}{\jumpdg{\v}}{\facedgint(D)}, 
    &&B_{p,D}^i(\u, \v)=\scaldg{\dgstab\jumpdg{\u}}{\jumpdg{\v}}{\facedgint(D)},\\
    &B_{c,D}^\partial(\u, \v)=\scaldg{\nu\nabla_h \u}{\v  \normal}{\facedgbdry(D)}, 
    &&B_{p,D}^\partial(\u, \v)=2\scaldg{\dgstab \u  \normal}{\v  \normal}{\facedgbdry(D)},
\end{align*}
where $\normal$ is the elementwise unit outward normal vector, $\scaldg{\cdot}{\cdot}{\facedgint(D)}$ and $\scaldg{\cdot}{\cdot}{\facedgbdry(D)}$ 
denote the $L^2$-inner products integrated over the union of all faces in $\facedgint(D)$ 
and $\facedgbdry(D)$, respectively.
Further, 
$
\dgstab=\frac{\gamma_0^2}{h_{\face}}  \frac{2\nu_1\nu_2}{\nu_1+\nu_2}
$
with a stabilization parameter, $\gamma_0$, and the face diameter, $h_{\face}$.
Note that $\jumpdg{\u}$ describes a jump term by definition.

If $D=\Omega$, we drop the subscript $D$ in the notations above.
The discrete formulation of (\ref{eq:PDEweak}) is to find $\uhe\in V_h$ such that 
$\ah(\uhe,\vh)  = \mathcal{F}(\vh)$ for all $\vh\in V_h$.
For later use, we define 
$D^+$ as the union of all $T\in\mesh$ with $T\cap D \neq \emptyset$ and $D^-$ as the union of all $T\in\mesh$ with $\overline{T}\cap D^c = \emptyset$, such that $D^- \subset D \subset D^+$. Throughout this paper, complements are taken with respect to $\Omega$.


\section{MS-GFEM for DG Discretizations of Elliptic PDEs}

Consider a mesh-resolved overlapping domain decomposition \(\{\omi\}_{\ddindex=1}^\numDom\), \(\cup_{\ddindex=1}^{\numDom} \omi = \Omega\), and define a partition of unity \(\{\chi_\ddindex\}_{\ddindex=1}^\numDom\) subordinate to this decomposition satisfying 
$\operatorname{supp}(\chi_\ddindex) \subset \overline{\omi^-}$, 
$0\leq \chi_{\ddindex}\leq 1$, 
$\sum_{\ddindex=1}^{\numDom}\chi_{\ddindex} \equiv 1$, 
$\chi_\ddindex\in\pcont{\Omega}{1}$. 
MS-GFEM builds local approximations on so-called oversampling domains $\osi$ satisfying $\omi\subset\osi\subset\Omega$ as the sum of a local source solution $u_\ddindex^p$ with $u_\ddindex^p=0$ on $\partial \osi \cap \partial \Omega$, cf.\ \cref{eq:particularFunction}, and a correction from an optimal $n_\ddindex$-dimensional spectral coarse space $S_{n_\ddindex}(\omi)$, cf.\ \cref{eq:locApproxSpace}. 
The global approximation is then assembled using a partition of unity,
$u^p := \sum_{\ddindex=1}^{\numDom} \chi_\ddindex u_\ddindex^p$, $S_{n}(\Omega) := \big\{\sum_{\ddindex=1}^{\numDom} \chi_\ddindex \phi_\ddindex : \phi_\ddindex \in S_{n_\ddindex}(\omi) \big\}$,
$n=\sum_{\ddindex=1}^\numDom n_\ddindex$,
and the MS-GFEM approximation is defined as $u^G = u^p + u^s$ with $u^s \in S_{n}(\Omega)$ satisfying 
$B(u^s,v) = \mathcal{F}(v) - B(u^p,v)$ for all $v \in S_{n}(\Omega)$. 
The key assumptions for the exponential error decay are a Caccioppoli inequality and a weak approximation property. We verify these, together with the remaining assumptions of \cite{Ma2025}, to show that our method fits into their general framework. For a detailed description of the abstract MS-GFEM and the complete statement of all assumptions, we refer the reader to \cite{Ma2025}.

\begin{proof}[Verification of {\cite[Assumption 2.3]{Ma2025}}]
(i) 
    For $D\subset\Omega$, we define 
    $
    \mathcal H(D)\coloneqq\{\vh|_D : \vh\in V_h\}
    $
    and
    $
    \mathcal H_0(D)\coloneqq\{\vh|_D : \vh\in V_h,  \vh = 0 \text{ on } 
    D\setminus D^-\},
    $
    both with the inner product
    \begin{align*}
        (\cdot,\cdot)_{\mathcal H(D)}
        = \scal{\nu\nabla_h\cdot}{\nabla_h\cdot}{D} 
        + \scaldg{\gamma_h^2\jumpdg{\cdot}}{\jumpdg{\cdot}}{\facedgint(D)}
        + \scaldg{\gamma_h^2\cdot}{\cdot}{\facedgbdry(D)}
        + \scal{\cdot}{\cdot}{D}.
    \end{align*}
    Clearly, $\mathcal H_0(D) \subset \mathcal H(D)$ for arbitrary subdomains $D\subset\Omega$, and $\mathcal H_0(\Omega)=\mathcal H(\Omega)$.\\
(ii)
    For $D\subset D^*$, we define 
    $
    E_{D,D^*}\colon\mathcal{H}_0(D)\to\mathcal{H}_0(D^*)
    $
    via $E_{D,D^*} (\vh) = \vh$ on $D^-$ and $E_{D,D^*} (\vh) = 0$ elsewhere. 
    We want to show that $\norm{E_{D,D^*} (\vh)}{\mathcal{H}_0(D^*)} = \norm{\vh}{\mathcal{H}_0(D)}$ 
    for all $\vh\in\mathcal{H}_0(D)$. Since $E_{D,D^*} (\vh)$ and $\vh$
    agree on $D^-$ and vanish elsewhere, the volume contributions to the 
    norms are equal. Next, we consider the contributions from interior 
    faces. 
    Let $F\in\facedgint$ with $F = \partial T_1 \cap \partial T_2$. We want to show 
    \begin{equation}
        \label{eq:extension_norm_equality_interior_faces}
    \scaldg{\gamma_h^2\jumpdg{E_{D,D^*}(\vh)}}{\jumpdg{E_{D,D^*}(\vh)}}{F\cap D^*}
        =  \scaldg{\gamma_h^2\jumpdg{\vh}}{\jumpdg{\vh}}{F \cap D}.
    \end{equation}
    If $F\subset D$, then $E_{D,D^*}(\vh)|_{T_1\cup T_2} = \vh|_{T_1\cup T_2}$ 
    and $F\cap D = F\cap D^*$, hence \cref{eq:extension_norm_equality_interior_faces} holds.
    If $F\subset D^c$, then $E_{D,D^*}(\vh)|_{T_1\cup T_2} = 0$
    and $F\cap D = \emptyset$, hence both sides of 
    \cref{eq:extension_norm_equality_interior_faces} vanish. 
    If $F\cap \partial D \neq \emptyset$, then $E_{D,D^*}(\vh)|_{T_1\cup T_2} =  \vh|_{T_1\cup T_2} = 0$
    and both sides of \cref{eq:extension_norm_equality_interior_faces} vanish.
    We proceed similarly for boundary faces $F\in\facedgbdry$ with $F \subset \partial T$ to show that
    \begin{equation}
        \label{eq:extension_norm_equality_boundary_faces}
        \scaldg{\gamma_h^2E_{D,D^*}(\vh)}{E_{D,D^*}(\vh)}{F\cap \partial D^*}
        =  \scaldg{\gamma_h^2\vh}{\vh}{F \cap \partial D}.
    \end{equation}
    If $T$ is not contained in $D^-$, then
    $E_{D,D^*}(\vh)|_T = 0$ and $\vh|_T = 0$ and thus 
    both sides of \cref{eq:extension_norm_equality_boundary_faces} vanish.
    If $T$ is contained in $D^-$,
    then $E_{D,D^*}(\vh)|_T = \vh|_T$ and since $F$ is a 
    boundary face with $T\subset D\subset D^*$, we have $F\subset \partial D$
    and $F\subset\partial D^*$. 
     Hence, $F\cap \partial D = F = F\cap \partial D^*$ and both sides of
    \cref{eq:extension_norm_equality_boundary_faces} are equal.\\
    (iii) 
        For $D\subset D^*$, we define 
        $R_{D^*,D}\colon \mathcal{H}(D^*) \to \mathcal{H}(D)$, $R_{D^*,D} (\vh) = \vh|_D$. 
        Due to the definition of the norm on $\mathcal{H}(D)$, we have
        $\norm{R_{D^*,D} (\vh)}{\mathcal H(D)} 
        \leq \norm{\vh}{\mathcal H(D^*)}$.\\
    (iv)
        For arbitrary $D\subset D^*$, $\uh \in \mathcal{H}(D^*)$ and $\vh \in \mathcal{H}_0(D)$, we have $B_D(\uh|_D,\vh) = B_{D^*}(\uh, E_{D,D^*}(\vh))$ 
        since $\vh$ vanishes on $D \setminus D^-$.
    \end{proof}
We denote by 
$
\interop\colon \{\vh : \vh|_T \in C^{\infty}(T) \text{ for all } T\in\mesh\} \to \mathcal{H}(\Omega)
$ 
the elementwise defined Lagrangian interpolation operator and define the partition of unity operator
$P_\ddindex\colon\mathcal{H}(\omi) \to \mathcal H_0(\omi)$, $\P_{\ddindex}(\u) = \interop(\chi_\ddindex \u)$. 
Boundedness of this operator is shown via the following lemma.

\begin{lemma}
    Let $\uh\in \mathcal{H}(\om)$, $\chi\in\pcont{\om}{1}$. Then, 
        $\norm{\chi\uh - \interop(\chi\uh)}{\mathcal{H}_0(\om)}
        \lesssim \norm{\chi\uh}{\mathcal{H}(\om)}$.
\end{lemma}

\begin{proof}
    First, we note that $\chi\uh$ is piecewise polynomial due to $\chi\in \pcont{\om}{1}$. 
    Therefore, an inverse inequality and elementwise interpolation properties of $\interop$ yield 
    \begin{align*}
        \lnorm{\nu^{1/2} \nabla_h(\chi\uh - \interop(\chi\uh))}{\om}
        &\lesssim \lnorm{\nu^{1/2} \nabla_h(\chi\uh)}{\om}, \\
        \lnorm{\chi\uh - \interop(\chi\uh)}{\om}
        &\lesssim  \lnorm{\chi\uh}{\om}.
    \end{align*}
    Next, we consider the jump terms. For all interior faces $F = \partial T_1 \cap \partial T_2$, observe that $\scaldg{\jumpdg{a}}{\jumpdg{b}}{\face} \leq (\lnorm{a|_{T_1}}{\face} + \lnorm{a|_{T_2}}{\face})(\lnorm{b|_{T_1}}{\face} + \lnorm{b|_{T_2}}{\face})$ for arbitrary $a$ and $b$, such that applying 
    the discrete trace inequality \cite[Lemma 1.46]{di2011mathematical}, elementwise interpolation properties and \cite[Lemma 1.43]{di2011mathematical} yields
    \begin{align}
            \notag&\scaldg{\dgstab\jumpdg{\left(\chi\uh- \interop(\chi\uh)\right)}}{\jumpdg{\left(\chi\uh- \interop(\chi\uh)\right)}}{\face}\\
            \notag&\lesssim \frac{\nu_1\nu_2}{\nu_1+\nu_2} h_{\tmin}^{-2} \lnorm{\chi\uh- \interop(\chi\uh)}{T_1 \cup T_2}^2 
            \lesssim \frac{\nu_1\nu_2}{\nu_1+\nu_2} h_{\tmin}^{-2} h_{\tmax}^2 \lnorm{\nabla_h(\chi\uh)}{T_1 \cup T_2}^2\\
            &\lesssim \frac{\nu_1\nu_2}{\nu_1+\nu_2} 
            \sum_{i=1}^{2}\frac{1}{\nu_i}\lnorm{\nu_i^{1/2}\nabla_h(\chi\uh)}{T_i}^2
             \leq \lnorm{\nu^{1/2}\nabla_h(\chi\uh)}{T_1 \cup T_2}^2, \label{eq:lemma1:face_estimate}
    \end{align}
    where $h_{\tmin} = \min\{h_{T_1},h_{T_2}\}$ and $h_{\tmax} = \max\{h_{T_1},h_{T_2}\}$.
    Summing over all interior faces of $\om$, we conclude
    \begin{align*}
        \scaldg{\gamma_h^2\jumpdg{(\chi\uh- \interop(\chi\uh))}}{\jumpdg{(\chi\uh- \interop(\chi\uh))}}{\facedgint(\om)}
        \lesssim \lnorm{\nu^{1/2}\nabla_h(\chi\uh)}{\om}^2.
    \end{align*}
    Boundary faces can be treated similarly.
    Finally, we combine all estimates to obtain the statement of the lemma.
\end{proof}

\begin{corollary}
    \label{cor:stability_interpol_new}
    Let $\uh\in \mathcal{H}(\om)$ and $\chi\in\pcont{\om}{1}$. Then,
    \begin{equation*}
        \norm{\interop(\chi\uh)}{\mathcal{H}_0(\om)}
        \lesssim \norm{\chi\uh}{\mathcal{H}(\om)}
        \lesssim \sqrt{1 + \norm{\nabla_h\chi}{L^{\infty}(\om)}^2} \norm{\uh}{\mathcal{H}(\om)}.
    \end{equation*}
\end{corollary}

\begin{proof}[Verification of {\cite[Assumption 2.9]{Ma2025}}] Let $\psi_\ddindex\in \mathcal H_0(\osi)$ be the unique solution of 
\begin{equation}\label{eq:particularFunction}
    \ahosj(\psi_\ddindex,\v) = \mathcal{F}_{\osi}(\v) \qquad \forall \v\in\mathcal H_0(\osi)
\end{equation}
and define the particular solution as $u_\ddindex^p \coloneqq \psi_\ddindex|_{\omi}$.
Note that $\ahosj$ does not contain face integrals over the interior subdomain boundary. However, since we define problem \cref{eq:particularFunction} on $\mathcal H_0(\osi)$, it is well-posed.
\end{proof}

\begin{proof}[Verification of {\cite[Assumption 2.13]{Ma2025}}] For every domain $\omi \subset D \subset \osi$, we define
$
    B_D^+(\cdot,\cdot) 
    \coloneqq (\nu\nabla_h\cdot, \nabla_h\cdot)_{D} + \scaldg{\gamma_h^2\jumpdg{\cdot}}{\jumpdg{\cdot}}{\facedgint(D)} 
    + \scaldg{\gamma_h^2\jumpdg{\cdot}}{\jumpdg{\cdot}}{\facedgbdry(D)}.
$
If $B_D^+(\uh,\vh) = 0$ holds for all $v\in\mathcal{H}(D)$, then 
this holds in particular for $\vh=\uh$. Thus, the gradient terms and jumps over inner faces in the definition of $B_D^+$ imply that $\uh$ is constant on $D$.
Conversely, if $\uh$ is constant, the jump terms and gradients in the definition of $B_D^+$ vanish and $B_D^+(\uh,\vh) = 0$ holds for all $\vh\in\mathcal{H}(D)$.
Hence, \cite[Assumption 2.13]{Ma2025} holds with $\mathcal{K}_D = \{0\}$ if $D$ is a boundary subdomain and with $\mathcal{K}_D = \mathrm{span}(1)$ if $D$ is an interior subdomain.
\end{proof}

For the construction of the coarse space, we define the space of discretely locally harmonic functions, 
$\mathcal{H}_B(\osi) \coloneqq \{\uh \in \mathcal{H}(\osi) : \ahosj(\uh,\vh) = 0 \; \forall \; \vh\in\mathcal{H}_0(\osi)\}$. Consider the generalized eigenvalue problem of finding 
$\eval \in [0,+\infty]$, $\evec \in \mathcal{H}_B(\osi)$ such that
$B^+_{\omi}(P_\ddindex(\evec|_{\omi}),P_\ddindex(\vh|_{\omi})) = \eval B^+_{\osi}(\evec,\vh)$ for all $\vh \in \mathcal{H}_B(\osi)$.
Denoting the $\evalindex$-th eigenpair as $(\evalloc, \evecloc)$, where $\evalloc[1] \geq \evalloc[2] \geq \ldots$, the local approximation space is built from the eigenfunctions corresponding to the $n_\ddindex$ largest eigenvalues:
\begin{equation}\label{eq:locApproxSpace}
  S_{n_\ddindex}(\omi) := \mathrm{span} \big\{ \evecloc[1]|_{\omi}, \ldots, \evecloc[n_\ddindex]|_{\omi} \big\}.
\end{equation}
Next, we verify the two central assumptions of \cite{Ma2025}. For $D\subset \Omega$, we
denote by $\Ccont(D)$ and $\Ccoer(D)$ the continuity and coercivity constants of the
bilinear form $B_{D}$ with respect to the (semi-)norm $\|\cdot\|_{B^+, D}$ induced by $B^+_D$.
Note that $\Ccont(D)$ and $\Ccoer(D)$ are independent of $\nu$, which can be shown as in 
\cite[Lemma 4.51, Lemma 4.52]{di2011mathematical}.

\begin{lemma}[Caccioppoli inequality, {\cite[Assumption 3.1]{Ma2025}}]
\label{lem:caccioppoli}
    Let $\om \subset \os \subset\Omega$ and $\u\in\mathcal{H}_B(\omega^{\ast})$ 
    with $\delta = \mathrm{dist}(\om, \partial \os\setminus\partial\Omega)>3 \max\limits_{K\in\mesh : K\cap \os\setminus \om \neq \emptyset} h_K$.
    Then, 
    \begin{equation}\label{eq:caccioppoli}
        \norm{\u|_{\om}}{B^+, \om}
        \lesssim \numax^{1/2} \delta^{-1} \lnorm{\u}{\os \setminus \om}. 
    \end{equation}
\end{lemma}

\begin{proof}
    Let $\eta\in\pcont{\omega^*}{1}$ be a cut-off function
    with $\operatorname{supp}(\eta) \subset \overline{(\os)^-}$, 
    $\eta = 1$ on $\om^+$ and $|\nabla_h\eta| \leq C_{\eta}\delta^{-1}$.  
    We proceed as in \cite{Ma2022error}. By definition,
    \begin{align}
        \label{eq:caccio_proof_1_new}
            \ahos(\eta\u,\eta\u) 
            = (\nu\nabla_h(\eta\u), \nabla_h(\eta\u))_{\os}
             + \ahosi(\eta\u,\eta\u)
            + \ahosb(\eta\u,\eta\u).
    \end{align}
    For the second term in \cref{eq:caccio_proof_1_new}, we use that $\eta$ is continuous 
    along mesh faces to compute
    \begin{align*}
        &\phantom{{}={}}\ahosi(\eta\u,\eta\u) 
        = \scaldg{\dgstab\jumpdg{\eta\u}}{\jumpdg{\eta\u}}{\facedgint(\os)}
        - \scaldg{\sumdgw{\nu\nabla_h(\eta\u)}}{\jumpdg{(\eta\u)}}{\facedgint(\os)} \\
        &= \scaldg{\dgstab\jumpdg{\u}}{\jumpdg{\eta^2\u}}{\facedgint(\os)}
        - \scaldg{\sumdgw{\nu\eta\nabla_h\u + \nu \u\nabla_h\eta}}{\jumpdg{(\eta\u)}}{\facedgint(\os)}\\ 
        &= \aposi(\u,\eta^2\u)
        - \acosi(\u,\eta^2 \u) 
        - \scaldg{\sumdgw{\frac{1}{2}\nu\eta\nabla_h\u + \nu\u\nabla_h\eta}}{\jumpdg{(\eta\u)}}{\facedgint(\os)}\\  
        &= \aposi(\u,\eta^2\u)
        - \acosi(\u,\eta^2 \u) 
        - \acosi(\eta^2 \u,\u)
        = \ahosi(\u,\eta^2\u).
    \end{align*}
    With the same arguments, we obtain
    $\ahosb(\eta\u,\eta\u) = \ahosb(\u,\eta^2\u)$. 
    Since $\eta$ is supported on $(\os)^-$, the harmonicity of $\u$ implies
    $\ahos(\u,\interop(\eta^2\u)) = 0$, 
    which we subtract from \cref{eq:caccio_proof_1_new} while plugging in the above identities and using the product rule 
    to obtain
    \begin{equation}
    \begin{aligned}
        \label{eq:caccio_proof_5_new}
        \ahos(\eta\u,\eta\u) 
            &= (\nu\u\nabla_h\eta, \u\nabla_h\eta)_{\os}
            +(\nu\nabla_h\u, \nabla_h(\eta^2\u - \interop(\eta^2\u)))_{\os} \\
            &\phantom{={}}+ \ahosi(\u,\eta^2\u - \interop(\eta^2\u)) 
            + \ahosb(\u,\eta^2\u - \interop(\eta^2\u)).
    \end{aligned}
    \end{equation}
    We will now show that 
    \begin{equation*}
        B_{\omega^*}(\eta\u,\eta\u) 
        \leq C \numax^{1/2} \lnorm{\u}{\os\setminus\om}
        \Big(\frac{1}{\delta} \lnorm{\nu^{1/2}\nabla_h(\eta\u)}{\os\setminus\om}+\frac{\numax^{1/2}}{\delta^2}\lnorm{\u}{\os\setminus\om}\Big).
    \end{equation*}
    Once we have this inequality, we can use a weighted Young's inequality to absorb the term involving the gradient of $\eta\u$ and infer
    \begin{equation}
        \label{eq:caccio_proof_young argument}
        B_{\omega^*}(\eta\u,\eta\u) 
        \leq \left(C + \frac{C^2}{2\Ccoer(\Omega)}\right) \frac{\numax}{\delta^2} \lnorm{\u}{\os\setminus\om}^2
        + \frac{\Ccoer(\Omega)}{2} \lnorm{\nu^{1/2}\nabla_h(\eta\u)}{\os\setminus\om}^2.
    \end{equation}
    Since $\eta$ is supported on $(\os)^-$, we have
    $B(\eta\u,\eta\u) = B_{\os}(\eta\u,\eta\u)$  
    and thus 
    $\Ccoer(\Omega) \norm{\eta\u}{B^+,\Omega}^2
        \leq B_{\os}(\eta\u,\eta\u)$ due to the coercivity of $B$. 
    Combining this with \cref{eq:caccio_proof_young argument}, the assumptions on $\eta$, and noting that the face integrals in the definition of $B^+_{\om}$ are bounded by the face integrals in the definition of $B^+_{\Omega}$ establishes the lemma.

    Now, we bound all terms in \cref{eq:caccio_proof_5_new}. 
    The first term can be bounded from above by
    $(\nu\u\nabla_h\eta, \u\nabla_h\eta)_{\os} 
    \leq \numax C_{\eta}^2 \delta^{-2}\lnorm{\u}{\os\setminus\om}^2$. 
    Since $\eta|_{\om^+} = 1$, we infer that all terms in \cref{eq:caccio_proof_5_new} only have contributions from elements contained in 
    $(\os)^-\setminus \om^+ \neq \emptyset$. Let $T\in\mesh$ be such an element. Using \cite[Lemma 5.4]{Ma2025} and an inverse inequality \cite[Lemma 1.44]{di2011mathematical}, 
    we obtain
    \begin{align*}
        &\phantom{{}\lesssim{}}(\nu\nabla_h\u, \nabla_h(\eta^2\u - \interop(\eta^2\u)))_{T}
        \leq \nu|_{T} \lnorm{\nabla_h\u}{T} \lnorm{\nabla_h(\eta^2\u - \interop(\eta^2\u))}{T}\\
        &\lesssim \nu|_{T} \lnorm{\nabla_h\u}{T}  \left( \frac{h_T}{\delta} \lnorm{\nabla_h(\eta\u)}{T}
        +\frac{h_T}{\delta^2}\lnorm{\u}{T} \right)\\
        &\lesssim \numax^{1/2}\lnorm{\u}{T}   \left( \frac{1}{\delta} \lnorm{\nu^{1/2}\nabla_h(\eta\u)}{T}
        +\frac{\numax^{1/2}}{\delta^2}\lnorm{\u}{T} \right)
    \end{align*}
    and hence 
    \begin{equation}\label{eq:caccio_proof_100_new}
    \begin{aligned}
        &\phantom{{}\lesssim{}}(\nu\nabla_h\u, \nabla_h(\eta^2\u - \interop(\eta^2\u)))_{\os}\\
        &\lesssim \numax^{1/2} \lnorm{\u}{\os\setminus\om}
        \left( \frac{1}{\delta} \lnorm{\nu^{1/2}\nabla_h(\eta\u)}{\os\setminus\om}+\frac{\numax^{1/2}}{\delta^2}\lnorm{\u}{\os\setminus\om} \right).
    \end{aligned}
    \end{equation}
    It remains to bound the last two terms of \cref{eq:caccio_proof_5_new}. 
    We only consider interior faces since the boundary terms can be treated in the same way.
    Recall that by definition
    $\ahosi(\u,\v)
                = \aposi(\u,\v)
                - \acosi(\u,\v)
                - \acosi(\v,\u)$.
    Consider two adjacent elements
    $T_1,T_2\in\mesh$ with common face $\face\in\facedgint$ and
    abbreviate $T\coloneqq T_1\cup T_2$ and $\nu_i=\nu|_{T_i}$. We bound
    the penalty term $\aposi$ face by face:
    \begin{align}
            \notag&\scaldg{\dgstab\jumpdg{\u}}{\jumpdg{(\eta^2\u- \interop(\eta^2\u))}}{\face}
            \lesssim \frac{\nu_1\nu_2}{\nu_1+\nu_2}  h_{\tmin}^{-2} \lnorm{\u}{T} \lnorm{\eta^2\u- \interop(\eta^2\u)}
            {T} \\
            \notag&\lesssim \frac{\nu_1\nu_2}{\nu_1+\nu_2} h_{\tmin}^{-2} \lnorm{\u}{T}
            \left( \frac{h_{\tmax}^2}{\delta} \lnorm{\nu^{1/2}\nabla_h(\eta\u)}{T}
            +\frac{h_{\tmax}^2}{\delta^2}\lnorm{\u}{T} \right)\\
            &\lesssim \numax^{1/2} \lnorm{\u}{T}
            \left( \frac{1}{\delta} \lnorm{\nu^{1/2}\nabla_h(\eta\u)}{T}
            +\frac{\numax^{1/2}}{\delta^2}\lnorm{\u}{T} \right),\label{eq:caccio_proof_123_new}
    \end{align}
    where the first step is carried out analogously to the first step in \cref{eq:lemma1:face_estimate}, and we further used 
    \cite[Lemma 5.4]{Ma2025} as well as \cite[Lemma 1.43]{di2011mathematical}.
    Note that the condition $\eta\in\pcont{\os}{1}$ is crucial for the applicability of the discrete trace inequality, 
    since it guarantees that $\eta^2\u$ is piecewise polynomial.
    Summing over all faces, we obtain the same upper bound for $\aposi(\u,\eta^2\u - \interop(\eta^2\u))$ as in \cref{eq:caccio_proof_100_new}. 

    For the first consistency term $\acosi(\u,\eta^2\u - \interop(\eta^2\u))$ of $\ahosi$, we again proceed as in \cref{eq:lemma1:face_estimate} to obtain
    \begin{align*}
            \scaldg{\sumdgw{\nu\nabla_h\u}}{\jumpdg{\eta^2\u - \interop(\eta^2\u)}}{\face}
            &\lesssim \frac{\nu_1 \nu_2}{\nu_1 + \nu_2}  h_{\tmin}^{-1} \lnorm{\nabla_h\u}{T}
            \lnorm{\eta^2\u - \interop(\eta^2\u)}{T}\\
            &\lesssim \frac{\nu_1 \nu_2}{\nu_1 + \nu_2}  h_{\tmin}^{-2} \lnorm{\u}{T}
            \lnorm{\eta^2\u - \interop(\eta^2\u)}{T},
    \end{align*}
    which can then be treated in the same way as \cref{eq:caccio_proof_123_new}, 
    such that summing over all faces once again yields the same upper bound for $\acosi(\u,\eta^2\u - \interop(\eta^2\u))$ as in \cref{eq:caccio_proof_100_new}. The term $\acosi(\eta^2\u - \interop(\eta^2\u),\u)$ can be treated analogously. 
\end{proof}

\begin{remark}
    Assuming that $\nu$ is piecewise constant is not strictly necessary but leads to 
    smaller constants in \cref{cor:stability_interpol_new,lem:caccioppoli}. 
    For general diffusion coefficients, an additional factor $\sqrt{\numax/\numin}$
    occurs in the estimates, cf. \cite[Lemma 3.10]{Ma2022error}.
\end{remark}

\begin{lemma}[Weak approximation property, {\cite[Assumption 3.4]{Ma2025}}]
    \label{lem:discrete_aharmonic_approximation_result_new}
    Let $\omega\subset \omega^{\ast}\subset {\omega}^{\ast\ast}$ be subdomains of $\Omega$ with 
    $\delta\coloneqq{\rm dist} ({\omega}^{\ast},\partial \omega^{\ast\ast}\setminus \partial \Omega)>\max\limits_{K\in\mesh : K\cap \os\setminus \om \neq \emptyset} h_K$, 
    and let $\mathsf{V}_{\delta}({\omega}^{\ast}\setminus \omega)\coloneqq \big\{{\bm x}\in {\omega}^{\ast\ast}: {\rm dist}({\bm x}, {\omega}^{\ast}\setminus \omega)\leq \delta \big\}$. 
    Then, there exists a constant $C_1>0$ depending only on $d$, such that
    for each integer $m\geq C_{1}\big|\mathsf{V}_{\delta}({\omega}^{\ast}\setminus \omega)\big|\delta^{-d}$ and 
    $h\leq \big|\mathsf{V}_{\delta}({\omega}^{\ast}\setminus \omega)\big|^{1/d} m^{-1/d}$, 
    there exists an $m$-dimensional space $Q_{m}(\omega^{\ast\ast})\subset L^{2}(\omega^{\ast\ast})$ such that 
    for all $\uh\in\mathcal{H}_B(\omega^{\ast\ast})$,
    \begin{align*}
        \inf_{\v\in Q_{m}(\omega^{\ast\ast})} \lnorm{\uh - \v}{{\omega}^{\ast}\setminus \omega}
        \lesssim \numin^{-1/2} \big|\mathsf{V}_{\delta}({\omega}^{\ast}\setminus \omega)\big|^{1/d} m^{-1/d} \norm{\uh}{B^+,\omega^{\ast\ast}}.
    \end{align*}
\end{lemma}

\begin{proof}
    By the assumption on the mesh size, we have 
    $\delta^-\coloneqq{\rm dist} ({\omega}^{\ast},\partial (\omega^{\ast\ast})^-\setminus \partial \Omega)>0$. 
    We define the set
    $\mathsf{V}_{\delta^-}({\omega}^{\ast}\setminus \omega)\coloneqq \big\{{\bm x}\in ({\omega}^{\ast\ast})^-: 
    {\rm dist}({\bm x}, {\omega}^{\ast}\setminus \omega)\leq \delta^- \big\}$  
    and denote by $\reconop\colon V_h((\omega^{\ast\ast})^-)\to \Winf{(\omega^{\ast\ast})^-}$ the reconstruction operator introduced in \cite[Section 3]{buffa2009compact} defined on $(\omega^{\ast\ast})^-$.
    By \cite[Lemma 5.5]{Ma2025}, there is an $m$-dimensional space
    $Y_m\subset L^{2}((\omega^{\ast\ast})^-)$ and a constant $C>0$ depending only on $d$ such that
    \begin{equation*}
        \inf_{\v\in Y_m}\Vert \reconop\uh-\v\Vert_{L^{2}(\omega^{\ast}\setminus \omega)} \leq C
        \big|\mathsf{V}_{\delta^-}({\omega}^{\ast}\setminus \omega)\big|^{1/d} m^{-1/d}\Vert \nabla_h \reconop\uh \Vert_{L^{2}((\omega^{\ast\ast})^-)}
    \end{equation*}
    for all $\uh\in \mathcal{H}((\omega^{\ast\ast})^-)$. Using the triangle inequality, 
    the previous estimate and \cite[Theorem 3.1]{buffa2009compact}, we obtain
    \begin{align*}
        &\inf_{\v\in Y_m} \lnorm{\uh - \v}{\omega^{\ast}\setminus \omega}
        \leq \lnorm{\reconop\uh - \uh}{\omega^{\ast}\setminus \omega}
        + \inf_{\v\in Y_m}  \lnorm{\reconop\uh - \v}{\omega^{\ast}\setminus \omega}\\
        &\lesssim h \norm{\uh}{B^+,(\omega^{\ast\ast})^-}
        + C_1
        \big|\mathsf{V}_{\delta^-}({\omega}^{\ast}\setminus \omega)\big|^{1/d} m^{-1/d}\Vert \nabla_h \reconop\uh \Vert_{L^{2}((\omega^{\ast\ast})^-)}\\
        &\lesssim \big|\mathsf{V}_{\delta^-}({\omega}^{\ast}\setminus \omega)\big|^{1/d} m^{-1/d} \numin^{-1/2}\norm{\uh}{B^+, (\omega^{\ast\ast})^-}.
    \end{align*}
    Defining $Q_m(\om^{\ast\ast})$ via extension by zero of functions in $Y_m$ and noting that 
    $\big|\mathsf{V}_{\delta^-}({\omega}^{\ast}\setminus \omega)\big| \leq \big|\mathsf{V}_{\delta}({\omega}^{\ast}\setminus \omega)\big|$, 
    the lemma is established.
\end{proof}

Having verified all relevant assumptions of \cite{Ma2025}, we obtain the desired nearly exponential decay of the eigenvalues and the global MS-GFEM error. Note that the Kolmogorov $n$-widths given in \cite[Theorems 3.8 and 2.23]{Ma2025} coincide with the eigenvalues stated below due to \cite[Lemmas 2.19 and 2.20]{Ma2025}.

\begin{theorem}[{\cite[Theorem 3.8]{Ma2025}}]
    \label{thm:exponential_decay}
    Assume that $\omi\subset\osi$ are truncated
    concentric cubes of side length $H_\ddindex$ and $H_\ddindex^*$, respectively, that are resolved by the mesh. Let $H_\ddindex^* > H_\ddindex$. Then, there exist constants $N_\ddindex\in\mathbb{N}$, $C_\ddindex > 0$ and $c_\ddindex > 0$, independent 
    of $h$, such that for all $n\geq N_\ddindex$, if $h$ is sufficiently small, we have
    \begin{equation*}
        \sqrt{\evalloc[n]}
        \leq C_\ddindex e^{-c_\ddindex n^{1/d}}.
    \end{equation*}
\end{theorem}

\begin{theorem} 
    The MS-GFEM solution $\uhG$ satisfies
    \begin{equation*}
        \norm{\uhe - \uhG}{B^+, \Omega}
        \leq \frac{\Ccont({\Omega})}{\Ccoer({\Omega})} \sqrt{\omcon\oscon} \left(\max_{\ddindex=1,\ldots,\numDom}\sqrt{\evalloc[n_{\ddindex+1}]}
        \frac{\Ccont(\osi)}{\Ccoer(\osi)}\right) \norm{\uhe}{B^+, \Omega},
    \end{equation*}
    where $\omcon$ and $\oscon$ are the coloring constants of \(\{\omi\}_{\ddindex=1}^\numDom\) and \(\{\osi\}_{\ddindex=1}^\numDom\), respectively.
\end{theorem}
\begin{proof}
    This follows from combining Céa's Lemma and \cite[Theorem 2.8]{Ma2025}.
    The latter is applicable with 
    $\varepsilon_\ddindex = \sqrt{\evalloc[n_{\ddindex+1}]} \Ccont(\osi) / \Ccoer(\osi) \norm{\uhe}{B^+, \osi}$ 
    due to \cite[Theorem 2.23]{Ma2025}, since
    $\norm{\uhe|_{\osi} - \psi_\ddindex}{B^+, \osi}
        \leq \Ccont(\osi) / \Ccoer(\osi) \norm{\uhe}{B^+, \osi}$
    and $\|\cdot\|_{B^+,\Omega}\leq \|\cdot\|_{\mathcal{H}(\Omega)}
    $.
\end{proof}

\printbibliography

@article{buffa2009compact,
  title={Compact embeddings of broken Sobolev spaces and applications},
  author={Buffa, Annalisa and Ortner, Christoph},
  journal={IMA Journal of Numerical Analysis},
  volume={29},
  number={4},
  pages={827--855},
  year={2009},
  publisher={OUP},
  DOI = {10.1093/imanum/drn038}
}

@Book{di2011mathematical,
  author    = {Di Pietro, Daniele Antonio and Ern, Alexandre},
  publisher = {Springer},
  title     = {Mathematical Aspects of Discontinuous Galerkin Methods},
  year      = {2011},
  ISBN = {9783642229800}
}

@article{Ma2022error,
  title = {Error estimates for discrete generalized FEMs with locally optimal spectral approximations},
  ISSN = {1088-6842},
  url = {http://dx.doi.org/10.1090/mcom/3755},
  DOI = {10.1090/mcom/3755},
  journal = {Mathematics of Computation},
  publisher = {American Mathematical Society (AMS)},
  author = {Ma,  Chupeng and Scheichl,  Robert},
  year = {2022},
  month = jul 
}

@article{Ma2025,
  title = {A Unified Framework for Multiscale Spectral Generalized FEMs and Low-Rank Approximations to Multiscale PDEs},
  ISSN = {1615-3383},
  url = {http://dx.doi.org/10.1007/s10208-025-09711-z},
  DOI = {10.1007/s10208-025-09711-z},
  journal = {Foundations of Computational Mathematics},
  publisher = {Springer Science and Business Media LLC},
  author = {Ma,  Chupeng},
  year = {2025},
  month = apr 
}
\end{document}

\typeout{get arXiv to do 4 passes: Label(s) may have changed. Rerun}